\documentclass[12pt]{amsart}
\usepackage{amsfonts}
\usepackage{amsmath, amscd, amssymb, yhmath}
\usepackage{graphicx,epsfig}
\usepackage[all,import,rotate]{xy}
\usepackage{latexsym}
\usepackage{mathrsfs}
\usepackage{xy}
\input xy
\usepackage{url}
\xyoption{all}

\newtheorem{theorem}{Theorem}[section]
\newtheorem{theorem1}{Theorem}
\newtheorem{lem}[theorem]{Lemma}

\newtheorem{cor}[theorem]{Corollary}
\newtheorem{prop}[theorem]{Proposition}
\newtheorem*{B-principle}{B-principle}
\newtheorem*{H-principle}{H-principle}

\theoremstyle{remark}
\newtheorem{remark}[theorem]{Remark}

\theoremstyle{definition}
\newtheorem{definition}[theorem]{Definition}

\begin{document}

\title{New link invariants and Polynomials (II), unoriented case }

\author{Zhiqing Yang, Jifu Xiao}
\thanks{The first author is supported by a grant (No. 10801021/a010402) of NSFC and grant
from Specialized Research Fund for the Doctoral program of Higher
Education (SRFDP) (20070141035)
}
\address{Department of Mathematics, Dalian University of Technology, China}
\email{yangzhq@dlut.edu.cn}

\date{\today}

\begin{abstract}
Given any unoriented link diagram, a group of new knot invariants are constructed. Each of them satisfies a generalized 4 term skein relation. The coefficients of each invariant is from a commutative ring. Homomorphisms and representations of such a ring defines new link invariants. In this sense, they produces the well-known Kauffman bracket, the Kauffman 2-variable polynomial, and the $Q$-polynomial. 
\end{abstract}


\maketitle

\setcounter{tocdepth}{2}
\tableofcontents


\section{Introduction}
Polynomial invariants of links are very well known. In 1928, J.W. Alexander {\cite{Al}} discovered the famous Alexander polynomial. 50 years later, in 1984 Vaughan Jones {\cite{J}} discovered the Jones polynomial. Soon, the HOMFLYPT polynomial was found {\cite{HOMFLY}\cite{PT}}. It turns out to be a generalization of both the Alexander polynomial and the Jones polynomial. Those polynomials satisfy some three term skein relations of the form $aL_++bL_-+cL_0=0$. There are other polynomials, for example, the Kauffman bracket, the Kauffman 2-variable polynomial, the $Q$-polynomial {\cite{K}}. Let's use the following symbols to denote the local link diagram patterns.

\begin{figure}[ht]
{\epsfxsize = 0.5 in $E_1$  {\epsfbox{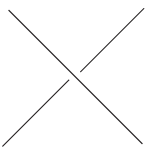}}
\epsfxsize = 0.5 in \quad \quad $E_2$ {\epsfbox{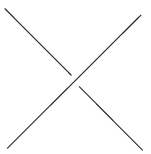}}
\epsfxsize = 0.5 in \quad \quad $E_0$ {\epsfbox{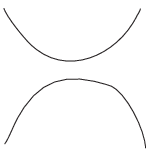}}
\epsfxsize = 0.5 in \quad \quad $E_{\infty}$ {\epsfbox{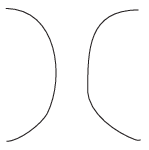}}
\caption{Local Diagrams with old notations}
\label{Figure 1}}
\end{figure}

\noindent Kauffman used the following equations to define link invariants.

\noindent (1) $E_1=aE_0+bE_{\infty}$

\noindent (2) $E_1+E_2=z(E_0+E_{\infty})$.

They are also semi-oriented invariants with the following properties. (1) The link invariant does not satisfy linear skein relation. (2) A writhe modification of some diagram invariant is the link invariant. (3) The diagram invariant is defined on unoriented diagrams, and satisfies a 4 term skein relation. A natural questions is, can they be further generalized?

In {\cite{Y}}, the first author defined new link invariants on oriented diagrams. In this paper, we will use similar method to construct new link invariants on unoriented link diagrams. It unifies the above semi-oriented link invariants. For simplicity, we switch to different system now.

\begin{figure}[ht]
{\epsfxsize = 3 in   {\epsfbox{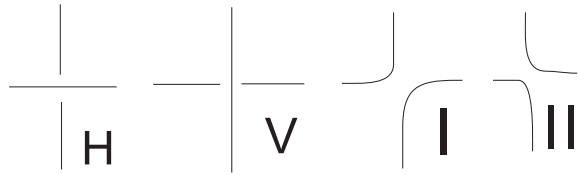}}
\caption{New notations in this paper}
\label{Figure 2}}
\end{figure}

Here, H denotes that the horizontal line is over the vertical, V denotes that the vertical line is over. I means type one smoothing, or smoothing in the first.

\noindent So the Kauffman equations can be written in a new form.

\noindent (1) $H=aI+bII$

\noindent (2) $H+V=z(I+II)$.

\noindent We propose the following skein relations:

\noindent If the two arrows/strands are from same link component, then
$$H+eV+aI+bII=0$$

\noindent If the two arrows/strands are not from same link component, then
$$H+e'V+a'I+b'II=0$$

\begin{remark} For simplicity, we use $E_1$ to denote both the link diagram with a special local pattern and the value of our invariant on the diagram $E_1$. But some times, when necessary, we use $f(E_1)$ to denote the value of our invariant on the diagram $E_1$.
\end{remark}

\noindent Let $B_1$ be the commutative ring generated by $e,e',a,a',b,b',v_n, n\geq 1$, and with the following relation set.

\noindent $\{e^2=1,e'{}^2=1,b=ea,b'=e'a', (e'-e)a'=(e'-e)a=0,aa=aa',(ee'-1)aa=0,(1+e+ea)v_n+av_{n+1}=0, n\geq 1\}$.

\bigskip $B_2$ has generator $a,a',b,b',v_n, n\geq 1$ and has the following relation set.

\noindent $\{a'b=ab'=ab,aa=aa',bb=bb',(1-a)a=(1-b)b,(1-a)v_n=bv_{n+1}, v_n=(a'+b')v_{n-1}, n\geq 1\}$.

Here are our main theorems.

\begin{theorem1} For unoriented link diagrams, there is a link invariant with values in $B_1$ and satisfies the following skein relations:

\noindent (1) If the two strands are from same link component, then
$H+eV+aI+bII=0$

\noindent (2) Otherwise,
$H+e'V+a'I+b'II=0$

The value for trivial n-component link is $v_n$. As in the oriented case, any homomorphism of $B_1$ also defines a link invariant.

\end{theorem1}

\begin{theorem1} For unoriented link diagrams, there is a link invariant $f(D)$ with values in $B_2$ and satisfies the following skein relations:

\noindent (1) If the two strands are from same link component, then
$H=aI+bII$

\noindent (2) Otherwise,
$H=a'I+b'II$

The value for a trivial n-component link diagram is $v_n$. As in the oriented case, any homomorphism of $B_2$ also defines a link invariant.
\end{theorem1}

There are also modified (by writhe) version of the two invariants. Like the oriented case, any homomorphic image of $B_1$ or $B_2$ also gives a new link invariant. For example, let $e=e'=1,a=a',b=b'$ in $B_1$, one can get the Kauffman 2-variable polynomial and the $Q$-polynomial. Similarly, $B_2$ produces the Kauffman bracket and the Jones polynomial. The rings $B_1,B_2,B_1',B_2'$ are easy to handle here, the word problems are solvable. 

\section{Type 1 invariant}

As mentioned before, we propose the following skein relations:

\noindent If the two arrows/strands are from same link component, then
$$H+eV+aI+bII=0$$

\noindent If the two arrows/strands are not from same link component, then
$$H+e'V+a'I+b'II=0$$

\begin{remark} In this paper, a link diagram invariant mean this invariant is well-defined for a fixed diagram, no matter which crossing point you resolve first, but may it may change under Reidemeister moves. A link invariant means a link diagram invariant which is also invariant under Reidemeister moves.
\end{remark}

\subsection{Resolution Consistence}

Given a local crossing, there is an ambiguity when apply the skein relation. The problem is that you can either regard it as $H$ or $V$ depending from which way you are looking at it. We have to solve this problem. The idea is to make the two results equal. From $H+eV+aI+bII=0$ you get $H=-\{eV+aI+bII\}$, but if you rotate the diagram 90 degree, then you get another equation $eH+V+bI+aII=0$. From this we shall get $H=-e^{-1}\{V+bI+aII\}$. Hence we ask $eV+aI+bII=e^{-1}\{V+bI+aII\}$ to be always true. The easiest solution is to ask $e=e^{-1}$, $a=e^{-1}b$, $b=e^{-1}a$. Those equations can be reduced to $e^2=1$, $ea=b$. If we multiply $H+eV+aI+bII=0$ by $e$ from left, we shall get $eH+V+bI+aII=0$, which is just the second equation above! The symmetry is clear here. Therefor, the equations $e^2=1$ and $ea=b$ guarantee that no matter which way we use the first skein  relation we shall get same result. Hence we don't use $eH+V+bI+aII=0$ in this paper.

\begin{remark}
There might be another choice. We ask $(e^2-1)V+(ea-b)I+(eb-a)II=0$ and $H+eV+aI+bII=0$ be both true, and use those two equations to define an invariant. But we have not figure it out yet.
\end{remark}

From now on, we use the following skein equations.

\noindent If the two arrows are from same component, then
$$H+eV+aI+bII=0$$

\noindent If the two arrows are not from same component, then
$$H+e'V+a'I+b'II=0.$$

We also ask $e'{}^2=1,e'a'=b'$.

\subsection{The algebra $B_1$}

Like in the oriented case, in order to get a link invariant, the symbols should satisfies certain relations. We work on the equations $f_{pq}=f_{qp}$ {\cite{Y}}.

Suppose we have two crossing points $p,q$. The easy case are that when one resole $p$, the crossing pattern (the information telling you which skein equation to use) at $q$ is not changed, and also when one resole $q$, the crossing pattern at $p$ is not changed. In this case, the results requires that any two elements of the following set commute.

$$\{e,a,b,e',a',b'\}$$

This relation set is denoted by $R^{B_1}_1$. The relations we discussed earlier: $e'{}^2=1,e'a'=b'$ and $e^2=1,ea=b$
will be denoted by $R^{B_1}_0$.

The nontrivial cases here are that the four ends of the two arcs of $p$ must be connected the two arcs of $q$. One can easily list the possibilities bellow:

\begin{figure}[ht]
{\epsfxsize = 4.5 in \centerline{\epsfbox{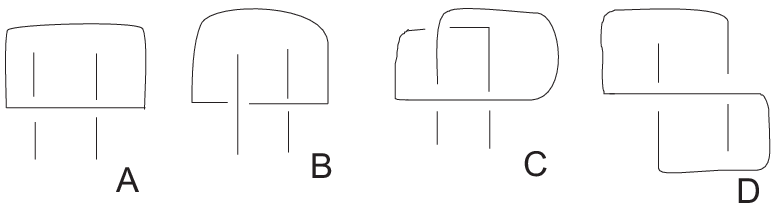}}}
\caption{}
\label{Figure 1}
\end{figure}

Here we do not draw the full picture. For example, case A actually represents two cases: $A_1,A_2$. They give same relations to the algebra.

\begin{figure}[ht]
{\epsfxsize = 3 in \centerline{\epsfbox{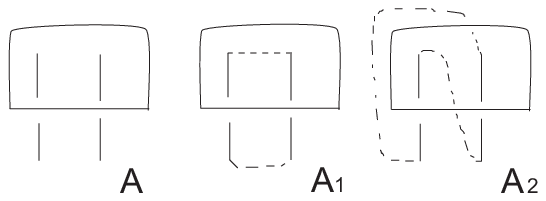}}}
\caption{}
\label{Figure 1}
\end{figure}

\noindent {\bf Case A}:

When we resolve the left crossing point first, we get
$$
\aligned (H,H)&=e'(V,H)+a'(I,H)+b'(II,H)) \\
 &= e'\{e'(V,V)+a'(V,I)+d'(V,II)\}\\
 &\quad\ +a'\{e(I,V)+a(I,I)+b(I,II)\} \\
 &\quad\ +b'\{e(II,V)+a(II,I)+b(II,II)\}
\endaligned
$$
When we resolve the left crossing point first, we get
$$
\aligned (H,H)&=e'(H,V)+a'(H,I)+d'(H,II)) \\
 &= e'\{e'(V,V)+a'(I,V)+b'(II,V)\}\\
 &\quad\ +a'\{e(V,I)+a(I,I)+b(II,I)\} \\
 &\quad\ +b'\{e(V,II)+a(I,II)+b(II,II)\}
\endaligned
$$
The two matrices are
\begin{tabular}[pos]{|l|l|l|l|l|l|l|l|l| }
\hline  1 $\backslash $ 2 & $V$ & $I$ & $II$      \\
\hline $V$ & $e'e'$ & $e'a'$ & $e'b'$     \\
\hline $I$ & $a'e$ & $a'a$ & $a'b$      \\
\hline $II$ & $b'e$ & $b'a$ & $b'b$   \\
\hline
 \end{tabular} \ \ \
\begin{tabular}[pos]{|l|l|l|l|l|l|l|l|l| }
\hline  1 $\backslash $ 2 & $V$ & $I$ & $II$      \\
\hline $V$ & $e'e'$ & $ea'$ & $eb'$     \\
\hline $I$ & $a'e'$ & $aa'$ & $cb'$      \\
\hline $II$ & $b'e'$ & $ba'$ & $bb'$   \\
\hline
 \end{tabular}

\noindent {\bf Case B}:

When we resolve the left crossing point first, we get
$$
\aligned (V,H)&=e'(H,H)+b'(I,H)+a'(II,H)) \\
 &= e'\{e'(H,V)+a'(H,I)+b'(H,II)\}\\
 &\quad\ +b'\{e(I,V)+a(I,I)+b(I,II)\} \\
 &\quad\ +a'\{e(II,V)+a(II,I)+b(II,II)\}
\endaligned
$$
When we resolve the left crossing point first, we get
$$
\aligned (V,H)&=e'(V,V)+a'(V,I)+b'(V,II)) \\
 &= e'\{e'(H,V)+b'(I,V)+a'(II,V)\}\\
 &\quad\ +a'\{e(H,I)+a(I,I)+a(II,I)\} \\
 &\quad\ +b'\{e(H,II)+a(I,II)+a(II,II)\}
\endaligned
$$
The two matrices are
\begin{tabular}[pos]{|l|l|l|l|l|l|l|l|l| }
\hline  1 $\backslash $ 2 & $V$ & $I$ & $II$      \\
\hline $H$ & $e'e'$ & $e'a'$ & $e'b'$     \\
\hline $I$ & $b'e$ & $b'a$ & $b'b$      \\
\hline $II$ & $a'e$ & $a'a$ & $a'b$   \\
\hline
 \end{tabular} \ \ \
\begin{tabular}[pos]{|l|l|l|l|l|l|l|l|l| }
\hline  1 $\backslash $ 2 & $V$ & $I$ & $II$      \\
\hline $H$ & $e'e'$ & $ea'$ & $eb'$     \\
\hline $I$ & $b'e'$ & $aa'$ & $bb'$      \\
\hline $II$ & $a'e'$ & $aa'$ & $ab'$   \\
\hline
 \end{tabular}

\noindent {\bf Case C}:

When we resolve the left crossing point first, we get
$$
\aligned (H,H)&=e(V,H)+a(I,H)+b(II,H)) \\
 &= e\{e(V,V)+a(V,I)+b(V,II)\}\\
 &\quad\ +a\{e(I,V)+a(I,I)+b(I,II)\} \\
 &\quad\ +b\{e'(II,V)+a'(II,I)+b'(II,II)\}
\endaligned
$$
When we resolve the left crossing point first, we get
$$
\aligned (H,H)&=e(H,V)+a(H,I)+b(H,II)) \\
 &= e\{e(V,V)+a(I,V)+d(II,V)\}\\
 &\quad\ +a\{e'(V,I)+a'(I,I)+b'(II,I)\} \\
 &\quad\ +b\{e(V,II)+a(I,II)+b(II,II)\}
\endaligned
$$
The two matrices are
\begin{tabular}[pos]{|l|l|l|l|l|l|l|l|l| }
\hline  1 $\backslash $ 2 & $V$ & $I$ & $II$      \\
\hline $V$ & $ee$ & $ea$ & $eb$     \\
\hline $I$ & $ae$ & $aa$ & $ab$      \\
\hline $II$ & $be'$ & $ba'$ & $bb'$   \\
\hline
 \end{tabular} \ \ \
\begin{tabular}[pos]{|l|l|l|l|l|l|l|l|l| }
\hline  1 $\backslash $ 2 & $V$ & $I$ & $II$      \\
\hline $V$ & $ee$ & $e'a$ & $eb$     \\
\hline $I$ & $ae$ & $a'a$ & $ab$      \\
\hline $II$ & $be$ & $b'a$ & $bb$   \\
\hline
 \end{tabular}

\noindent {\bf Case D}:

When we resolve the left crossing point first, we get
$$
\aligned (H,H)&=e(V,H)+a(I,H)+d(II,H)) \\
 &= e\{e(V,V)+a(V,I)+b(V,II)\}\\
 &\quad\ +a\{e'(I,V)+a'(I,I)+b'(I,II)\} \\
 &\quad\ +b\{e(II,V)+a(II,I)+b(II,II)\}
\endaligned
$$
When we resolve the left crossing point first, we get
$$
\aligned (H,H)&=e(H,V)+a(H,I)+b(H,II)) \\
 &= e\{e(V,V)+a(I,V)+b(II,V)\}\\
 &\quad\ +a\{e'(V,I)+a'(I,I)+b'(II,I)\} \\
 &\quad\ +b\{e(V,II)+a(I,II)+b(II,II)\}
\endaligned
$$
The two matrices are
\begin{tabular}[pos]{|l|l|l|l|l|l|l|l|l| }
\hline  1 $\backslash $ 2 & $V$ & $I$ & $II$      \\
\hline $V$ & $ee$ & $ea$ & $eb$     \\
\hline $I$ & $ae'$ & $aa'$ & $ab'$      \\
\hline $II$ & $be$ & $ba$ & $bb$   \\
\hline
 \end{tabular} \ \ \
\begin{tabular}[pos]{|l|l|l|l|l|l|l|l|l| }
\hline  1 $\backslash $ 2 & $V$ & $I$ & $II$      \\
\hline $V$ & $ee$ & $e'a$ & $eb$     \\
\hline $I$ & $ae$ & $a'a$ & $ab$      \\
\hline $II$ & $be$ & $b'a$ & $bb$   \\
\hline
 \end{tabular}

\bigskip\noindent Compare those matrices, we get the following relation set :

\noindent
$R^{B_1}_2=\{(e'-e)a'=(e'-e)b'=(e'-e)a=(e'-e)b=0,a'b=ab'=ab,aa=aa',bb=bb'\}$

\subsection{Construction and Proves}

To prove the above setting produces a link invariant, we need oriented the link diagrams. Later, we shall prove orientation independence.

Given any link diagram $D$, we shall first assume add the following information.

   (1) Suppose each link component has an {\bf orientation}.

   (2) Give an {\bf order} of the components by integers: 1,2, $\cdots$, m.

   (3) On each component $k_i$, pick a {\bf base point} $p_i$.

Now, we go through component $k_1$ from $p_1$ along its orientation. When we finish $k_1$, we shall pass to $k_2$ start from $p_2$, $\cdots$.

\begin{definition} A crossing is called {\bf bad} if it is first passed over, otherwise, it is called {\bf good}. A link diagram contains only good crossings is called a {\bf monotone/ascending} diagram.
\end{definition}

\begin{lem} A monotone diagram can be monotonously reduced to a zero crossing diagram using Reidemeister moves without increasing crossing number at each step.
\end{lem}

\begin{cor} A monotone diagram corresponds to a trivial link.
\end{cor}

\bigskip We shall define the invariant inductively on the index pair $(c,d)$, where $c$ is the crossing number of the diagram, and $d$ is the number of bad points of the diagram. It is obvious that $d\leq c$.

\begin{prop} The invariant satisfies the following properties.

\noindent(0) The value for any link diagram is uniquely defined.

\noindent(1) Satisfying skein relations if we resolve at any bad point.

\noindent(2) Invariant under base point change.

\noindent(3) Invariant under Reidemeister moves.

\noindent(4) Invariant under changing order of components.

\noindent(5) Invariant under orientation change.

\end{prop}

\bigskip\noindent {\bf Proof of the statement (0)}: As in paper 1, we do inductions on index $(c,d)$, where $c$ is crossing number, $b$ is the number of bad points. $0\leq d\leq c$.

\bigskip\noindent {\bf Step 1}. For a diagram of index $(n,0)$, define its value to be $v_n$.

\noindent Then the claim (0)-(5) is satisfied for diagram with index $(n,0)$.

Now suppose the claim (0)-(5) is prove for link diagrams with index strictly less than $(c,d)$.

\noindent {\bf Step 2}. If the diagram $D$ has bad points, say its index is $(c,d)$, where $d>0$, we resolve the diagram at the first bad point $p$. Then, in the skein equation, all the other terms are of smaller indices than $(c,d)$. By induction hypothesis, the terms with smaller $c$ can be defined and are invariant under base point change, orientation change, and changing order of components. So we can choose the base points, orientation, and ordering of link components arbitrarily. There is one term corresponds to crossing change, and it has a canonical base point set and ordering of link components. So all the terms except one in the skein equation have been uniquely defined, hence the skein relation uniquely defines the value for $D$.

$\\$\noindent {\bf Proof of the claim (1)}:

For a link diagram $D$, if $D$ has at most one bad point, then by definition, it satisfies claim (1).
If $D$ has at least 2 bad points, and one resolve at a bad point $q$. If $q$ is the first bad point, then by definition, the equation is satisfied.
If not, denote the first bad point by $p$. Denote the value of $D$ by $f(D)$. If we resolve at $p$, we get many diagrams $D_1,D_2,\cdots $. Let $f_p(D)$ denote the signed weighted sum of those diagrams. Then by definition $f(D)=f_p(D)$. Each diagram $D_i$, has lower indices than $(c,b)$. We resolve each $D_i$ at $q$, then we get the  signed weighted sum  $f_{q}(D_i)$. By induction hypothesis, $f(D)=f_p(D)=\sum f_{q}(D_i)$.

On the other hand, we can resolve $D$ at $q$ first, we get many diagrams $D_1',D_2',\cdots $, each has lower indices than $(c,b)$. Hence the claim (1)-(4) are satisfied. We get a signed weighted sum $f_q(D)$. We resolve each $D_i'$ at $p$, then we get the  signed weighted sum  $f_p(D_i')$. By induction hypothesis, $f_q(D)=\sum f_p(D_i')$. However, the algebra is designed such that $\sum f_p(D_i')=\sum f_{q}(D_i)$! (This is the equation $F_{pq}=F_{qp}$.)

Therefor, $f(D)=f_p(D)=\sum f_{q}(D_i)=\sum f_p(D_i')=f_q(D)$. That is, if we resolve at $q$, the skein equation is satisfied.

\begin{cor} Resolving a link diagram at any point (not necessarily bad point), the skein equation is satisfied.
\end{cor}

\begin{proof} If $q$ is a good point of $D$, we make a crossing change at $q$ get a new diagram $D'$, then $q$ is bad point of $D'$. If we resolve $D'$ at $q$,  the skein equation is satisfied. But this the same equation of $D$ resoling at $q$.
\end{proof}

\bigskip\noindent {\bf Proof of the claim (2)}:

Given a diagram $D$ with a fixed ordering of components and orientation, suppose that there are two base point set $B$ and $B'$. We only need to deal with the case that $B$ and $B'$ has only one point $b$ and $b'$ different, they are in the same component $k$, and between $b$ and $b'$ there is only one crossing point $p$. In the base point systems $B$ and $B'$, $D$ has the same bad points except $p$. If there is bad point other than $p$, say $q$, we resolve $D$ at $q$ to get diagrams $D_1,D_2,\cdots $. Then those $D_i$'s has lower indices than $D$, hence base point invariance is proved for them. On the other hand, the skein equation is proved, hence before resolving, the values for $D$ with different base point systems are the same.

If there is no other bad points, there are two cases. Case 1. $p$ is a good point for both the two base point systems, then the values for $D$ are both $v_n$, hence equal.
Case 2. $p$ is a bad point for both the two base point systems, then the skein equation tells the values are the same.

Case 3, $p$ is good in $B$, bad in $B'$.
Then The diagram in $B$ is a monotone diagram. Use a similar argument as in lemma 1 one can prove that we can fix the crossing $p$ and monotonously reduce all other crossings by Reidemeister moves( denote those moves by $\Omega$). The proof is similar, and one use an outmost argument if necessary.

It follows that all the smoothings at $p$ produces trivial links no matter before or after the Reidemeister moves $\Omega$.
Then before the Reidemeister moves, in $B$, the value of $D$ is $v_n$.  In $B'$, the value is uniquely defined by the skein equation. If the value in $B'$ is $v_n$, then plug this into the skein equation, we have $v_n+ev_n+av_n+bv_{n+1}=0$.

On the other hand, this is also a sufficient condition. So, as long as the symbols always satisfy the equation $v_n+ev_n+av_n+bv_{n+1}=0$ for any $n\geq 1$, the value for $D$ in $B'$ is $v_n$, hence we proved base point invariance.

\bigskip\noindent {\bf Proof of the statement (3)}:

\noindent {\bf (i)} Given two diagrams $D$ and $D'$, which differs at a Reidemeister move I. Say $D$ has index $(c,d)$, where $D'$ has index $(c+1,d')$. In the local Reidemeister move I part, say the crossing point is $p$. $D$ and $D'$ have the same bad points except $p$. Like in (2),  if there is bad point other than $p$, we can resolve and prove Reidemeister move I invariance inductively.

Otherwise, the other points are all good, then $D$ and $D'$ are both diagrams of trivial links. If $p$ is good in $D'$, there is nothing to prove. If
$p$ is bad, we can use (2), base point invariance, to get rid of this bad point then get the proof.

$\\$\noindent {\bf (ii)} We shall prove invariance under Reidemeister move III first. Given two diagrams $D$ and $D'$, which differs at a Reidemeister move III.  Likewise, we can assume all other points are good. In the two local diagrams containing the Reidemeister move III, there is a one to one correspondence between the three arcs appearing in the two local diagrams. We can also order the three arcs by 1,2,3,($1',2,3'$ in $D'$) such that arc 1 is above arc, and arc 2 is above arc 3. The one to one correspondence preserves the ordering. Their intersections induce a one to one correspondence between the three pair points in the two diagrams.  Call them $p,p'$, $q,q',r,r'$.

Suppose $p$ is the intersection of arc 1 and arc 2 (or arc 2 and arc 3), then we can resolve $p,p'$. There are many terms of lower crossing. It is easy to see that those terms equal each other in pairs (Reidemeister move II with crossing number $c-1$ is proved). There is only one term left, which corresponds to crossing change. Therefor, if change the crossing $p,p'$ get two diagrams $D_p, D_{p'}'$, the value $D=D'$ if and only if value $D_p= D_{p'}'$.

This technique can't be applied to the intersection of arc 1 and arc 3, say $r$, directly. However, if $r$ is a bad point, then one can easily see that $p,q$ can't be both good. Say $p$ is bad. After resolving at $p$, one can resolve $r$.

Hence one can make all the three intersections good. It follows that $p',q,r'$ are good now. Now the invariance is clear.

\bigskip\noindent {\bf (iii)} Given two diagrams $D$ and $D'$, which differs at a Reidemeister move II. $D'$ has two more crossings $p,q$. Likewise, we can assume all other points are good. In the local picture contains the Reidemeister move II, if one crossing is good, the other is bad, one can use a base point change to make those two points good. Then both the diagrams are diagrams for trivial link. There is nothing to prove.

If both the two crossing are bad, and base point changes wouldn't change them from bad to good, changing both the two crossing will make them both good. Hence both the diagrams are diagrams for trivial link.

Now let's first work on the case that one diagram $D$ is a trivial, the other $D'$ has only two bad crossings. The crossings are intersections from two link components. We have the following diagrams.
\begin{figure}[ht]
{\epsfxsize = 3.5 in \centerline{\epsfbox{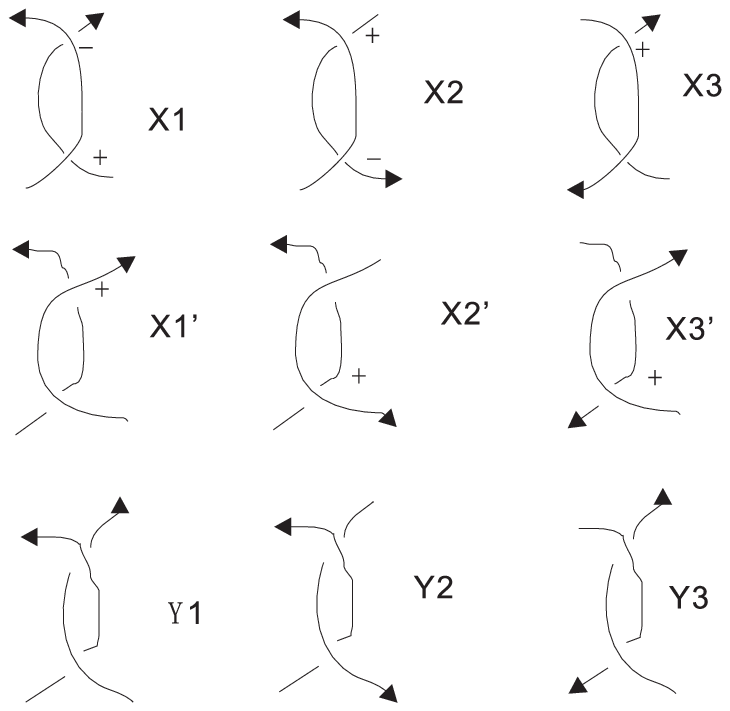}}}
\caption{}
\label{Figure 4}
\end{figure}

We shall show the diagrams $Xi$ and $Xi'$ have same value for $i=1,2,3$.
Let's resolve both $Xi$ and $Xi'$ at the positive crossing point, then we have

$X_i+eY_i+av_{n-1}+bv_{n-1}=0$ and $X_i'+eY_i+av_{n-1}+bv_{n-1}=0$. Hence we have that $Xi$ and $Xi'$ have same value for $i=1,2,3$.

By changing from $X_i$ to $X_i'$ or from $X_i'$ to $X_i$, we can get a monotone diagram. Hence we have $f(X_i)=f(X_i')=v_n=f(D)$.

In the general case, lemma shows that in all the diagrams $Xi,Xi',Yi$, we can monotonously reduce all other crossing points and at the same time keep the local diagram fixed. Hence the last 4 terms are diagrams of the trivial link of $n-1$ components. On the other hand, the last 4 terms in the skein equation has $c-1$ crossings, hence their Reidemeister invariance is proved. Therefor, their values are all $v_{n-1}$. Hence we also have $f(X_i)=f(X_i')=v_n=f(D)$.

\bigskip\noindent {\bf Proof of the statement (4)}: The induction is a little different here. Since the above proves do not need (4) for crossing $=c$, only need (4) for crossings $<c$ (to define the value on the smaller crossing diagrams), we assume finished proving Reidemeister moves invariance for all diagrams with crossing $< c+1$.

Given two diagrams $D$ with different ordering of components. For simplicity, call them $D^1$ and $D^2$ with indices $(c,d_1)$ and $(c,d_2)$. Suppose $d_1\leq d_2$. We inducts on $(c,d_1)$. If $d_1=0$, $D^1$ is a trivial link diagram, so is $D^2$. However, $D^2$ has bad points. Since we already proved Reidemeister move invariance. Using lemma 1, we can monotonously change the diagrams $D^1, D^2$ to trivial link diagrams of disjoint circles on the plane using Reidemeister moves. Hence $D^2$ also has value $v_n$.

If $d_1\neq 0$, then we resolve at a bad point. Now the invariance follows by induction.

\bigskip\noindent {\bf Proof of the statement (5)} Now calculation of any oriented diagram is reduce a linear combination of monotone diagram $D_1,D_2,\cdots $. This linear relation is orientation independent since the skein relation does not depend on orientation. The value of a monotone diagram is $v_n$ even you change its orientation since we have Reidemeister move invariance. This proves (5).

\subsection{Modify by writhe}

There is another closely related link invariant. The idea is that the skein relations can reduce the calculation to monotone diagrams, and we can regard the set of monotone diagrams as a basis and assign writhe dependant values on those diagrams. Now the skein relations don't give one a link invariant, but we can make a new function $g(w)$, such that the product $g(w)f(D)$ is Reidemeister invariant. Here $w$ is a the writhe of the link diagram. This is very similar to the Kauffman 2-variable polynomial.

\begin{prop} For oriented link diagrams, there are invariants satisfies the following properties.

\noindent (1) The value for monotone diagram $D$ is defined as $f(D)=h(w)v_n$, where $w$ is a the writhe of the link diagram, $n$ is the number of components.

\noindent (2) The value $f(D)$ for any link diagram $D$ is uniquely defined.

\noindent (3) $f(D)$ satisfies skein relations if we resolve at any bad point.

\noindent (4) $f(D)$ is invariant under base point change.

\noindent (5) There is another function $g(w)$ such that $g(w)h(w)=1$ and  $F(D)=g(w)f(D)$ is invariant under all Reidemeister moves. $F(D)$ is also invariant under base point change, and $f(D)$ is invariant under Reidemeister moves II and III.

\noindent (6) $F(D)$, and $f(D)$ are invariant under changing order of components.

\end{prop}

As before, the proof is an induction on index $(c,b)$. For statement (1), there is nothing to prove.

\bigskip\noindent {\bf Proof of the statement (2)(3)}: The proof is almost the same as last section. We also always resolve at the first bad point.

\bigskip\noindent {\bf Proof of the statement (4)}: Like before, we use induction on $(c,b)$, and we need an extra equation here. Suppose we have a diagrams with different base point sets $B,B'$. For simplicity, the diagram will be called $D$ with base point set $B$, and $D'$ with base point set $B'$. We only need to deal with the case that $B$ and $B'$ has only one point $x$ and $x'$ different, they are in the same component $k$, and between $x$ and $x'$ there is only one crossing point $p$. As before, we can assume there is no bad point except $p$. In the base point systems $B'$, $D'$ has one bad point $p$, and $D$ is a monotone diagram.

Hence $f(D)=h(w)v_n$. When we use skein relation to calculate $f(D')$, we get $f(D')+h(w-2)ev_n+ah(w-1)v_{n+1}+bh(w-1)v_n=0$. Hence we need the equation $h(w)v_n+h(w-2)ev_n+ah(w-1)v_{n+1}+bh(w-1)v_n=0$. this equation is sufficient to prove the base point invariance.

\bigskip\noindent {\bf Proof of the statement (5)}:

\noindent (i) For Reidemeister move I, the case is different now. We can still first reduce the proof to monotone diagrams $D,D'$, they have writhe $w,w+1$. $D$ has value $h(w)v_n$, $D'$ has value $h(w+1)v_n$. So we require $g(w)$ has the property that $g(w)h(w)=g(w+1)h(w+1)$. So we require $g(w)h(w)\equiv 1$, and $g(w),h(w)$ commute with all other symbols.

\noindent (ii) The proves for Reidemeister move III invariance of $F(D)$ is the same as before.

\noindent (iii) For Reidemeister move II invariance of $F(D)$, the result is a little different here. As before, we have two diagrams $D,D'$. $D$ is a monotone diagram, with $c$ crossing, and writhe is $w$. $D'$ has $c+2$ crossing. Those two bad crossings, say $p,q$ are intersections of different components. Then use the same argument as last proof for Reidemeister move II invariance, we get two equations. Last time, we had $X_i+e'Y_i+(a'+b')v_{n-1}=0$ and $X_i'+e'Y_i+(a'+b')v_{n-1}=0$. Now they should be modified a little bit, we have to add the writhe part into the equations. It is clear that the last 2 terms in each of the two equations have writhe $w-1$. So it is also true that $Xi$ and $Xi'$ have same value for $i=1,2,3$. So we can change $X_i$ to $X_i'$ or vise versa. After changing both $p,q$ to good points, the proof of Reidemeister move II invariance is trivial for $F(D)$.

\noindent (iv) Since Reidemeister move III and II does not change writhe, $F(D)=F(D')$ implies $f(D)=f(D')$. Hence $f(D)$ is invariant under Reidemeister move II and III.

\bigskip\noindent {\bf Proof of the statement (6)} Proof for $F(D)$ is also the same as before. Since $g(w)f(D)=F(D)$ is invariant under ordering change, so is $f(D)$.

\begin{remark}
The difference here is that although $F(D)=f(D)$ on monotone diagrams, $F(D)$ does not satisfy the skein relations.
\end{remark}
\begin{remark}
\noindent An easy choice for the equation $f(D')+h(w-2)ev_n+ah(w-1)v_{n+1}+bh(w-1)v_n=0$ is to let $h(w)=A^w$ for a new variable $A$. Then the equation is reduced to $A^2v_n+ev_n+aAv_{n+1}+bAv_n=0$.
\end{remark}

\section{Type 2 invariant}

   There is another new knot invariant.

   Same component, they satisfies the following relation:

    $H=aI+bII$

     \noindent Different components, they satisfies the following relation:

    $H=a'I+b'II$

  This invariant is related to the Kauffman bracket. There is another set of equations for them.

\subsection{The algebra $B_2$}

As in last section, we first require that any two elements from $\{a,b,a',b'\}$ commute. Then we study the nontrivial cases. It turns out we need to work on the same case as last section too. So we don't need to carry out the calculation again. The skein relation $H+eV+aI+bII=0$ can be reduced to $H=aI+bII$ by drop the $eV$ term and then change signs.

So the relations for this algebra is, in $R^{B_1}_2=\{(e'-e)a'=(e'-e)b'=(e'-e)a=(e'-e)b=0,a'b=ab'=ab,aa=aa',bb=bb'\}$, delete all the terms containing $e$ or $e'$, hence we get $R^{B_2}_2=\{a'b=ab'=ab,aa=aa',bb=bb'\}$

\subsection{Construction and proves}

\begin{prop} For unoriented link diagrams, there is an invariant satisfies the following properties.

\noindent (0) The value is defined uniquely for any unoriented link diagram.

\noindent (1) Satisfying skein relations if we resolve at any crossing point.

\noindent (2) Invariant under Reidemeister moves for any two diagrams with crossing $<c$.

\end{prop}

\bigskip\noindent {\bf Proof of the statement (0)}:
We shall define the invariant inductively on crossing number $c$ of the diagram.

\bigskip\noindent {\bf Step 1}. For a n-component oriented link diagram of crossing number $c=0$, define its value to be $v_n$.

\noindent Then the claim (0)-(2) is satisfied.

\noindent {\bf step 2}. If the diagram $D$ has crossing points, we resolve the diagram at one crossing point $p$. Then, in the skein equation, all the other terms are of smaller crossing numbers. By induction hypothesis, the other terms are uniquely defined. Hence the skein relation uniquely defines the value for $D$.

$\\$
\noindent {\bf Proof of the statement (1)}:

For a link diagram $D$, if $D$ has only one crossing point, we get the equations: $v_n=bv_n+av_{n+1}$ and $v_n=av_n+bv_{n+1}$.

To make those two equations consistent, we have introduce new relations on $a,b$. They can changed to $(1-a)v_n=bv_{n+1}$ and $(1-b)v_n=av_{n+1}$, so we ask $(1-a)a=(1-b)b$. This is also sufficient for consistency. The only problem is when one want to rewrite $adv_{n+1}$ in $v_n$, there are two ways. The two results are $b(1-b)v_n$ and $a(1-a)v_n$. Hence We need to add $(1-a)a=(1-b)b$ to our algebra relations.

If $D$ has at least two crossing points. The proof is the same as the proves before. We use one point $p$ to define the value of the diagram $D$. For any other crossing point $q$, we have to prove the skein equation.

Denote the value of $D$ by $f(D)$. If we resolve at $p$, we get many diagrams $D_1,D_2,\cdots $. Let $f_p(D)$ denote the signed weighted sum of those diagrams. Then by definition $f(D)=f_p(D)$. Each diagram $D_i$, has lower indices than $c$. We resolve each $D_i$ at $q$, then we get the  signed weighted sum  $f_{q}(D_i)$. By induction hypothesis, $f(D)=f_p(D)=\sum f_{q}(D_i)$.

On the other hand, we can resolve $D$ at $q$ first, we get many diagrams $D_1',D_2',\cdots $, each has lower indices than $c$. Hence the claim (1)-(2) are satisfied. We get a signed weighted sum $f_q(D)$. We resolve each $D_i'$ at $p$, then we get the  signed weighted sum  $f_p(D_i')$. By induction hypothesis, $f_q(D)=\sum f_p(D_i')$. However, the algebra is designed such that $\sum f_p(D_i')=\sum f_{q}(D_i)$! (This is the equation $F_{pq}=F_{qp}$.)

Therefor, $f(D)=f_p(D)=\sum f_{q}(D_i)=\sum f_p(D_i')=f_q(D)$. That is, if we resolve at $q$, the skein equation is satisfied.

$\\$\noindent {\bf Proof of the statement (2)}:

\noindent {\bf (i)} Given two diagrams $D$ and $D'$, which differs at a Reidemeister move I. Say $D$ has index $c$, where $D'$ has index $c+1$. In the local Reidemeister move I part, say the crossing point is $p$. $D$ and $D'$ have the same crossing points except $p$. If there are crossing points other than $p$, we can resolve both the diagrams and prove Reidemeister move I invariance inductively.

Otherwise, $p$ is the only crossing point, then $D$ and $D'$ are both diagrams of trivial links. Then\ Reidemeister move I invariance is guaranteed by the following equations
$(1-a)a=(1-b)b$, as in proof of statement 1.

$\\$\noindent {\bf (ii)} Given two diagrams $D$ and $D'$, which differs at a Reidemeister move II. Likewise, we can assume there is no other crossing points. There are two cases. Case 1. The Reidemeister move II involves only one link component. Then the invariance follows from Reidemeister move I invariance. Case 2. The Reidemeister move II involves two link components. Resolve one crossing point and use Reidemeister move I invariance one get $f(D)=(a'+b')v_{n-1}$. Then Reidemeister move II invariance follows from the following equation:
$v_n=(a'+b')v_{n-1}$.

$\\$\noindent {\bf (iii)} Given two diagrams $D$ and $D'$, which differs at a Reidemeister move III.  Likewise, we can assume all other points are good. In the local diagram containing the Reidemeister move III, there is a one to one correspondence between the three arcs appearing in the two local diagrams. We can also order the three arcs by 1,2,3,($1',2,3'$ in $D'$) such that arc 1 is above arc, and arc 2 is above arc 3. The one to one correspondence preserve the ordering. Suppose arc 1 and arc 2 intersects at $p$, arc $1'$ and arc $2'$ intersects at $p'$. Then we can resolve at $p,p'$ at the same time. The resulting terms can be paired up and equal each other since we prove Reidemeister move II invariance. Therefor, we proved Reidemeister move III invariance.

\subsection{Modify it by writhe}

As in type 1 invariant, the type 2 invariant can also be modified by writhe. Instead of asking the above definition to be Reidemeister moves invariant, we can ask its modification to be Reidemeister moves invariant. Denote the value of diagram $D$ by $f(D)$, let $\omega$ denote the writhe, $c$ denote the crossing number, $\mu$ denote number of link components. We ask a new family of functions $g(\omega , c, \mu)$ with parameters in  $\omega , c, \mu$, such that $gf(D)$ is Reidemeister moves invariant.

\begin{prop} There is an invariant satisfies the following properties.

\noindent (0) The value $f(D)$ is defined uniquely for any unoriented link diagram, and on trivial link diagrams ($c=0$) has value $v_n$.

\noindent (1) $f(D)$  satisfies skein relations if we resolve at any crossing point.

\noindent (2) $F(D)=g(\omega , c, \mu)f(D)$ is invariant under Reidemeister moves for any two diagrams with crossing $<c$.

\end{prop}

\bigskip\noindent {\bf Proof of the statement (0)(1)} Same as above.

\bigskip\noindent {\bf (i)} To make it Reidemeister move I invariant, like above, we need some new equations. Suppose there are two diagrams $D$ and $D'$, which differs at a Reidemeister move I. Say $D$ has index $c$, where $D'$ has index $c+1$. In the local Reidemeister move I part, say the crossing point is $p$. $D$ and $D'$ have the same crossing points except $p$. Suppose diagram $D$ has parameters $\omega , c, \mu$, the $D'$ has parameters $\omega +1 , c+1, \mu$ or $\omega -1, c+1, \mu$.

Like before, we can resolve all other crossings, and the resulting terms for $D$ and $D'$ can be paired up. Now we can resolve $p$. For example, if $p$ has positive crossing, then we can group the terms for $D'$ together, such that each group has the form $av_{n+1}+bv_n$ for some $n$, and for $D$, there is one term $v_n$ corresponds to it.

Hence we can add the following equations for Reidemeister move I invariance:

$g(\omega +1 , c+1, \mu)\{av_{n+1}+bv_n\}=g(\omega , c, \mu)v_n$.

Likewise, if $p$ has positive crossing, we get another equation.

$g(\omega -1 , c+1, \mu)\{bv_{n+1}+av_n\}=g(\omega , c, \mu)v_n$.

\bigskip\noindent {\bf (ii)} To make it Reidemeister move I invariant, like above, we need some new equations. Suppose there are two diagrams $D$ and $D'$, which differs at a Reidemeister move II. Say $D$ has index $c$, where $D'$ has index $c+2$. In the local Reidemeister move II part, say the crossing point is $p,q$. $D$ and $D'$ have the same crossing points except $p,q$. Suppose diagram $D$ has parameters $\omega , c, \mu$, the $D'$ has parameters $\omega , c+2, \mu$.

We can resolve all other crossings, and the resulting terms for $D$ and $D'$ can be paired up. Say $D_i$ and $D_i'$ is one of the pairs. $D_i$ has no crossings. We resolve $D_i'$ at the negative crossing point first, then the positive point, we always get $a'(av_{n}+bv_{n-1} ) +b'(a'v_{n} +b'v_{n-1})$

Then we can group the terms for $D'$ together, such that each group has the form $a'(av_{n-1}+bv_n) +b'(a'v_{n} +b'v_{n-1})$
 for some $n$, and for $D$, there is one term $v_n$ corresponds to it.

Hence we can add the following equations for Reidemeister move II invariance:

\noindent $g(\omega , c+2, \mu)\{a'(av_{n-1}+bv_n)+b'(a'v_{n}+ b'v_{n-1}) \}=g(\omega , c, \mu)v_n$.

\bigskip\noindent {\bf (iii)} For the Reidemeister move III invariance, things are much easier. The only difference of $D$ and $D'$ is position of a crossing point. The equality is trivial.

$\\$
In Kauffman's bracket, if we modify it by writhe, we get the Jones polynomial. This is the same idea here.

\section{Type 3 and type 4 invariant}

\noindent {\bf Type 3  invariant}

\noindent If the two arrows are from same component, then
$$H+eV+aI+bII=0$$

\noindent If the two arrows are not from same component, then
$$H=a'I+b'II$$

As before, we ask $b=ea$.

The second relation set for $B_3$ is a modification of
$R^{B_1}_2=\{(e'-e)a'=(e'-e)b'=(e'-e)a=(e'-e)b=0,a'b=ab'=ab,aa=aa',bb=bb'\}$. We delete all the terms has $e'$, so we have $R^{B_3}_2=\{ea'=eb'=ea=eb=0,a'b=ab'=ab,aa=aa',bb=bb'\}$. However, from $ea=eb=0$ and $ea=b$, we shall have $a=eb=0=ea=b$. The we have $H+eV=0$. This invariant is not of much interest.

$\\$
\noindent {\bf Type 4 invariant}

\noindent If the two arrows are from same component, then
$$H=aI+bII$$

\noindent If the two arrows are not from same component, then
$$H+e'V+a'I+b'II=0$$

As before, we ask $b'=e'a'$.

The second relation set for $B_3$ is a modification of
$R^{B_1}_2=\{(e'-e)a'=(e'-e)b'=(e'-e)a=(e'-e)b=0,a'b=ab'=ab,aa=aa',bb=bb'\}$. We delete all the terms has $e$, so we have $R^{B_3}_2=\{e'a'=e'b'=e'a=e'b=0,a'b=ab'=ab,aa=aa',bb=bb'\}$. However, from $e'a'=e'b'=0$ and $e'a'=b'$, we shall have $a'=e'b'=0=e'a'=b'$. The we have $H+e'V=0$. This invariant is not of much interest.

\section{Conclusion}

We have four algebras $B_1,B_1',B_2,B_2'$.

\noindent $B_1$ has generator $e,e',a,a',b,b',v_n, n\geq 1$ and with relation sets

\noindent $R^{B_1}_0=\{e^2=1,e'{}^2=1,b=ea,b'=e'a'\}$,

\noindent $R^{B_1}_1:$ any two elements of $e,e',a,a'$ commute.

\noindent $R^{B_1}_2=\{(e'-e)a'=(e'-e)a=0,aa=aa',(ee'-1)aa=0\}$,

\noindent $R^{B_1}_3=\{(1+e+ea)v_n+av_{n+1}=0, n\geq 1\}$.

\begin{theorem} For unoriented link diagrams, there is a link invariant with values in $B_1$ and satisfies the following skein relations:

\noindent (1) If the two strands are from same link component, then
$H+eV+aI+bII=0$

\noindent (2) Otherwise,
$H+e'V+a'I+b'II=0$

The value for trivial n-component link is $v_n$.

\end{theorem}

\begin{remark}
If one let $v_1=1$, $e=e'=1, a=a'$, then one get the $Q$-polynomial.
\end{remark}

\bigskip There are functions $h(w)$, $g(w)$ defined for $w\in Z$, such that $g(w)h(w)=1$.

\noindent $B_1'$ has generator $e,e',a,a',b,b',v_n, n\geq 1$ and with relation sets

\noindent $R^{B_1}_0=\{e^2=1,e'{}^2=1,b=ea,b'=e'a'\}$

\noindent $R^{B_1}_1:$ any two elements of $e,e',a,a',h(w),g(w)$ commute.

\noindent $R^{B_1}_2=\{(e'-e)a'=(e'-e)a=0,aa=aa',(ee'-1)aa=0\}$,

\noindent $R^{B_1}_3=\{[h(w)+h(w-2)e+h(w-1)ea]v_n+ah(w-1)v_{n+1}=0, n\geq 1\}$.

\begin{theorem} For unoriented link diagrams, there is a link diagram invariant with values in $B_1'$ and satisfies the following skein relations:

\noindent (1) If the two strands are from same link component, then
$H+eV+aI+bII=0$

\noindent (2) Otherwise,
$H+e'V+a'I+b'II=0$

The value for a monotone n-component link diagram is $h(w)v_n$, where $w$ is the writhe of the diagram. $F(D)=g(w)f(D)$ is a link invariant.
\end{theorem}

\begin{remark}
If one let $v_1=1,e=e'=-1$, $a=a',b=b'$, and  $h(w)=a^w, g(w)=a^{-w}$, then one get the Kauffman 2-variable polynomial.
\end{remark}

\bigskip $B_2$ has generator $a,a',b,b',v_n, n\geq 1$ and with relation sets

\noindent $R^{B_1}_0=\emptyset $

\noindent $R^{B_1}_1:$ any two elements of $a,a',b,b'$ commute.

\noindent $R^{B_1}_2=\{a'b=ab'=ab,aa=aa',bb=bb'\}$,

\noindent $R^{B_1}_3=\{(1-a)a=(1-b)b,(1-a)v_n=bv_{n+1}, v_n=(a'+b')v_{n-1}, n\geq 1\}$.

\begin{theorem} For unoriented link diagrams, there is a link invariant $f(D)$ with values in $B_2$ and satisfies the following skein relations:

\noindent (1) If the two strands are from same link component, then
$H=aI+bII$

\noindent (2) Otherwise,
$H=a'I+b'II$

The value for trivial n-component link $D$ is $v_n$.
\end{theorem}

\begin{remark}
Even let $v_1=1$, $a=a',b=b'$, one don't get the Kauffman Bracket. This is a new link invariant.
\end{remark}

\bigskip For type $B_2'$, there are functions $h(w,c,\mu)$, and $g(w,c,\mu)$ defined for $w\in Z,c,\mu \in N,$, and $g(w,c,\mu)h(w,c,\mu)=1$. It has generator $e,e',a,a',b,b',v_n, n\geq 1$ and with relation sets

\noindent $R^{B_1}_0=\emptyset$

\noindent $R^{B_1}_1:$ any two elements of $a,a',b,b',g(w,c,\mu)$ commute.

\noindent $R^{B_1}_2=\{a'b=ab'=ab,aa=aa',bb=bb'\}$,

\noindent $R^{B_1}_3=\{g(\omega +1 , c+1, \mu)\{av_{n+1}+bv_n\}=g(\omega , c, \mu)v_n,g(\omega -1 , c+1, \mu)\{bv_{n+1}+av_n\}=g(\omega , c, \mu)v_n, g(\omega , c+2, \mu)\{a'(av_{n-1}+bv_n)+b'(a'v_{n}+ b'v_{n-1}) \}=g(\omega , c, \mu)v_n,  n\geq 1\}$.

\begin{theorem} For unoriented link diagrams, there is a link diagram invariant $f(D)$ with values in $B_2'$ and satisfies the following skein relations:

\noindent (1) If the two strands are from same link component, then
$H=aI+bII$

\noindent (2) Otherwise,
$H=a'I+b'II$

$f(D)$ for a trivial n-component link diagram is $v_n$. $F(D)=g(\omega , c, \mu)f(D)$ is a link invariant.
\end{theorem}

\begin{remark}
This invariant can produce the Jones polynomial.
\end{remark}

Let's see why. If we let $a=a',b=b'$, and $g(\omega , c, \mu)=A^{\omega}$, the relations are reduce to the followings.

\noindent $R^{B_1}_1:$ any two elements of $a,b,A$ commute.

\noindent $R^{B_1}_3=\{A(av_{n+1}+bv_n)=v_n,bv_{n+1}+av_n=Av_n, (a^2+b^2)v_{n-1}+2abv_n=v_n,  n\geq 1\}$. If we ask $b=a^{-1}$, $A=-a^{-3}$, and $a=t^{-\frac{1}{4}}$, then we get the Jones polynomial. However, it seems unlikely this system of equations has another easy to handle solution.

The above invariants do not produce other interesting polynomials in the usual sense, for if one ask $a,b$ to have inverses, then he shall get $a=a',b=b'$. This is not like the oriented cases in  {\cite{Y}}. However, the generators and relations are much less here, hence we can directly work on the invariants themselves. For example, in $B_1$, we have the following relations.

$b=ea,b'=e'a'$, $e'a'=ea',e'a=ea, aa'=aa$, and $av_{n+1}=-(1+e+ea)v_n$.

We replace every word in the left hand side by the right hand side, we shall get a unique ``simplest'' representative. In this sense, the wrod problem is solvable here. Hence those invariants can be easily used directly.


\begin{thebibliography}{Ag}

\bibitem{A}  Colin C. Adams, {\it The Knot Book}, W.H. Freeman and Company (1999)

\bibitem {Al}  J.W. Alexander, {\it Topological Invariants of Knots and Links}, Transactions of the American Mathematical Society, Volume 30, Issue 2 (April 1928), 275-306.

\bibitem{GH} G. Burde, and H. Zieschang, {\it Knots}, de Gruyter Studies in Mathematics, 5, Walter de Gruyter, Berlin (1985).

\bibitem{HOMFLY} P. Freyd, D. Yetter, J. Hoste, W. B. R. Lickorish, K. Millett, and A. Ocneanu, {\it A new polynomial invariant of knots and links}, Bull. Amer. Math. Soc. 12 (1985), 239-249.

\bibitem{J}  V. F. R. Jones, {\it Hecke algebra representations of braid groups and link polynomials}, Ann. of Math. (2) 126 (1987), 335-388.

\bibitem{K}  L. Kauffman, {\it State models and the Jones polynomial}, Topology 26 (1987), 395-407.

\bibitem{M}  V. O. Manturov, {\it Knot Theory}. CRC Press, 2004.

\bibitem{PT} J. H. Przytycki and P. Traczyk, {\it Invariants of links of Conway type}, Kobe J. Math. 4 (1987), 115-139.

\bibitem{R} Rolfsen, D, {\it Knots and links}, Publish or Perish Inc., Berkeley, (1976), 160-197.

\bibitem{Y} Yang, Z, {\it New link invariants and Polynomials (I), oriented case}, preprint.

\end{thebibliography}
\end{document}